\documentclass[12pt]{article}

\usepackage[]{graphics,graphicx,float,latexsym,subfigure}
\usepackage{amsfonts,amstext,amsmath,amssymb,amsthm,times}

\newcommand{\ignore}[1]{} 

\usepackage[lmargin=1in,rmargin=1in,tmargin=1in,bmargin=1in]{geometry}
\numberwithin{equation}{section} 

\newtheorem{theorem}{Theorem}
\newtheorem{lemma}{Lemma}

\theoremstyle{definition} 
\newtheorem*{definition}{Definition}
\newtheorem*{remark}{Remark}

\newcommand{\mrm}[1]{\mathrm{#1}}
\newcommand{\Exp}{{\mrm{Exponential}}}
\newcommand{\Geom}{{\mrm{Geometric}}}
\newcommand{\ARL}{{\mrm{ARL}}}

\newcommand{\SD}{{\mrm{SD}}}

\newcommand{\mc}[1]{\mathcal{#1}} 
\newcommand{\Fc}{\mc{F}}

\newcommand{\mb}[1]{\mathbf{#1}}
\newcommand{\Pb}{{\mb{P}}}
\newcommand{\Eb}{{\mb{E}}}

\def\One{\mathchoice{\rm 1\mskip-4.2mu l}{\rm 1\mskip-4.2mu l}%
{\rm 1\mskip-4.6mu l}{\rm 1\mskip-5.2mu l}}
\newcommand\Ind[1]{{\One_{\{#1\}}}}

\newcommand{\La}{\Lambda}
\newcommand{\de}{\delta}

\newcommand{\ga}{\gamma}
\newcommand{\vae}{\varepsilon}

\newcommand{\tX}{\tilde{X}}
\newcommand{\hX}{\hat{X}}
\newcommand{\hN}{\hat{N}}

\newcommand{\abs}[1]{\left\vert#1\right\vert}
\newcommand{\set}[1]{\left\{#1\right\}}

\newcommand{\brc}[1]{\left(#1\right)}

\newcommand{\xra}{\xrightarrow}

\usepackage{titlesec} 

\titleformat{\section}{\large\bfseries}{\thesection.}{.5em}{}
\titlespacing{\section}{0pt}{*3}{*2}
\titleformat{\subsection}{\normalfont\bfseries}{\thesubsection.}{.5em}{}
\titlespacing{\subsection} {0pt}{7pt}{*1}
\titleformat{\subsubsection}[runin]{\normalfont\it\bfseries}{\thesubsubsection.}{.5em}{}[.]
\titlespacing{\subsubsection} {0pt}{5pt}{*1}

\usepackage{fancyhdr}

\begin{document}

\title{\large \bf ASYMPTOTIC EXPONENTIALITY OF THE DISTRIBUTION OF FIRST EXIT TIMES
FOR A CLASS OF MARKOV PROCESSES WITH APPLICATIONS TO QUICKEST CHANGE DETECTION}

\author{\small {\sc Moshe Pollak}\\
\small The Hebrew University of Jerusalem\\
\small Department of Statistics\\
\small Mount Scopus\\
\small Jerusalem 91905, Israel\\
\small msmp@mscc.huji.ac.il\\
 \and
\small {\sc Alexander G.\ Tartakovsky}\\
\small University of Southern California\\
\small Department of Mathematics \\
\small 3620 S.\ Vermont Ave, KAP-108\\
\small Los Angeles, CA 90089-2532, USA\\
\small tartakov@usc.edu}

\date{\normalsize \sl Submitted to Probability Theory and Its Applications, March 2007}

\maketitle

\begin{abstract}
We consider the first exit time of a nonnegative Harris-recurrent Markov process
from the interval $[0,A]$ as $A\to\infty$. We provide an alternative method of
proof of asymptotic exponentiality of the first exit time (suitably standardized)
that does not rely on embedding in a regeneration process. We show that under
certain conditions the moment generating function of a suitably standardized
version of the first exit time converges to that of $\Exp(1)$, and we connect
between the standardizing constant and the quasi-stationary distribution (assuming
it exists).  The results are applied to the evaluation of a distribution of run
length to false alarm in change-point detection problems.

\

\end{abstract}

\markboth{M. Pollak and A. G. Tartakovsky}{Asymptotic Exponentiality of First Exit
Times with Applications to Change Detection}

\section{Introduction} \label{s:Intro}

Let $(\Omega, \Fc, \Pb )$ be a probability space and $\{X(n)\}$, $n=0,1,2,\dots$ be
a discrete-time non-negative Harris-recurrent Markov process defined on this space.
The limiting distribution as $A\to \infty$ of the suitably standardized first exit
time of the process from the interval $[0,A]$ turns out often to be exponential.

The standard method for proving this asymptotic exponentiality is to try to find a
version of the process that is regenerative (cf.\ Glasserman and Kou, 1995 and
Asmussen, 2003). The heuristic behind this is that since the process is
Harris-recurrent, it returns to a given set over and over again, and thus creates
``cycles" that are ``almost independent." Hence, the first cycle in which $X(n)$
exceeds $A$ is approximately geometrically distributed, and if the expected length
of a cycle is finite and the probability of exceeding $A$ in a given cycle tends to
$0$ as $A\to\infty$, then, suitably standardized,  the asymptotic distribution of
the first exit time is exponential.

In this paper, we make a connection between the standardization constant and the
quasi-stationary distribution.  Our method of proof is a coupling argument.
Although less general as a method for proving asymptotic exponentiality than the
regeneration argument, we believe that our method is of interest in its own right.
This notwithstanding, the regeneration argument seems to be widely unknown in the
statistics community, and ought to be publicized.

The paper is organized as follows. In Section~\ref{s:Results}, we present the main
result that states that the limiting distribution of the suitably standardized
version of the first exit time as $A\to\infty$ is $\Exp(1)$ and that the moment
generating function converges to that of $\Exp(1)$, which implies that the
convergence is in $L^p$ for all $p \ge 1$. The proof is given in
Section~\ref{s:Proof}. We make a few remarks in Section~\ref{s:Remarks}. In
Section~\ref{s:Applications},  we give examples and describe applications to the
evaluation of the distribution of the run length to false alarm for several change
detection procedures.

\section{Main Results} \label{s:Results}

Let $\{X(n)\}_{n=0}^\infty$ be a discrete-time Harris-recurrent Markov process with state space
$[0,\infty)$ and stationary transition probabilities. Let $\Pb^x$ denote the probability measure
for the process when it starts at $x$ (i.e., $X(0)=x$), and let $\Pb^G$ denote the probability
measure when the initial state is distributed according to the distribution $G$.

\begin{definition}
We call the process stochastically monotone if $\Pb^x(X(1) \ge y)$ is non-decreasing
and right-continuous in $x$ for all $y$.
\end{definition}

We will be interested in the behavior of the first exit time of $X(n)$ from the interval $[0,A]$
when $X(n)$ starts at $x\in [0,A)$, i.e., of the stopping time
\begin{equation} \label{ST}
N_A^x= \min\set{n \ge 1: X(n) > A}, \quad X(0)=x,
\end{equation}
where $0\le x < A$ and $A$ is a positive finite threshold, assuming that the Markov process $X(n)$
is stochastically monotone and Harris-recurrent.

The following theorem is the main result of the paper.

\begin{theorem} \label{th1}
Let $X(n)$, $n=0,1,2,\dots$ be a stochastically monotone Harris-recurrent Markov process with state
space $[0,\infty)$ and stationary transition probabilities such that:

\textbf{C1.} The stationary distribution $H(y)=\lim_{n\to\infty} \Pb^x\set{X(n) \le
y}$ exists and its support is $[0,\infty)$.

\textbf{C2.} The quasi-stationary distribution $H_A(y)$ $=$ $\lim_{n\to\infty}
\Pb^x\set{X(n) \le y | N_A^x > n}$ exists for all $0 \le x < A$ and for all
$0<A<\infty$.

Let $p_A=\Pb^{H_A}\set{X(1) >A}$.

Then:

\textbf{(i)} The distribution of $p_A \, N_A^x$ is asymptotically $\Exp(1)$ as
$A\to \infty$ for all fixed $x \in [0,\infty)$.

\textbf{(ii)} The moment generating function $\Eb \exp\set{t p_A N_A^x}$ of $p_A
N_A^x$ converges to $1/(1-t)$ as $A\to\infty$ for all fixed $x \in [0,\infty)$. In
particular, it follows that
\[
\lim_{A\to\infty} p_A \Eb N_A^x =1 \quad \text{and} \quad
\lim_{A\to\infty}{\mrm{Variance}}\set{p_A N_A^x} =1 .
\]

\end{theorem}

Conditions C1 and C2 hold in a variety of scenarios. See corresponding remarks in Section
\ref{s:Remarks} and examples in Section \ref{s:Applications}.

We begin with a heuristic argument. A formal proof requires several auxiliary results and is given
in Section \ref{s:Proof}.

Write $N_A^{H_A}$ for the stopping time when the process $X(n)$ starts at a random
point $X(0)=\xi$ in $[0,A]$ that has a quasi-stationary distribution $H_A$, i.e.,
$\Pb(\xi \le y)= H_A(y)$. Then $\Pb^{H_A}(X(n) > A | N_A^{H_A} \ge n) = p_A$ for
all $n \ge 1$, and, therefore, the distribution of $N_A^{H_A}$ is geometric with
the parameter $p_A$ for all $A>0$. Further, under conditions C1 and C2, the
probability $p_A$ goes to 0 as $A\to\infty$, which implies that $p_A N_A^{H_A}$
converges weakly to $\Exp(1)$ as $A\to\infty$. Intuitively, the asymptotic behavior
of the stopping time $N_A^x$ for every fixed point $x$ is similar to that of
$N_A^{H_A}$. Mathematical details are presented in the next section.

\section{Proof} \label{s:Proof}

In order to prove Theorem \ref{th1} we need the following lemmas. We use the notation of the
previous section, and we assume that the conditions of Theorem \ref{th1} are satisfied.

\begin{lemma} \label{le1} 
The quasi-stationary distribution
\[
H_A(y)=\lim_{n\to\infty} \Pb^x\set{X(n) < y | N_A^x > n}
\]
converges to the stationary distribution $H(y)$ at all continuity points $y$ of $H$.
\end{lemma}

\begin{proof}
Follows directly from Theorem 1 of Pollak and Siegmund (1986).
\end{proof}


Recall that $N_A^{H_A}$ is the stopping time \eqref{ST} when the Markov process $X(n)$ starts from
the random point that has the quasi-stationary distribution $H_A$, i.e., $X(0)\sim H_A$.

\begin{lemma} \label{le2} 
 The distribution of $N_A^{H_A}$ is $\Geom(p_A)$, where $p_A=$  $\Pb^{H_A}\set{X(1) >A}$. Hence $p_A \Eb
N_A^{H_A}=1$ and $p_A N_A^{H_A}$ converges in distribution to $\Exp(1)$ as $A\to\infty$.
\end{lemma}

\begin{proof}
Since the Markov process is Harris-recurrent, there is no absorbing state, so that $\Pb(N_A^{H_A} =
\infty)=0$. Therefore, the geometric property of $N_A^{H_A}$ is obvious. Lemma \ref{le1} and the
assumption that the support of $H$ is $[0,\infty)$ guarantee that $p_A \xra[A\to\infty]{} 0$.
\end{proof}


\begin{lemma} \label{le3} 
Let $X^x(n)$ denote a process that starts from $x$ and has the same transition probabilities as
$X(n)$. Let $0 \le x < y < \infty$. There exists a sample space with $X^x(n)$ and $X^y(n)$ such
that $X^y(n) \ge X^x(n)$ for all $n \ge 1$.
\end{lemma}

\begin{proof}
Clearly $X^y(1)$ is stochastically larger than $X^x(1)$, so that one can construct a
sample space where $X^y(1) \ge X^x(1)$. To complete the proof, continue by induction
on $n$.
\end{proof}


\begin{lemma} \label{le4} 
Let $0 \le x < y < \infty$. Let $\tX^x(n)$ and $\tX^y(n)$ be independent Markov processes started
at $x$ and $y$ respectively, both having the same transition probabilities as $X(n)$. Then
\begin{equation} \label{prob1}
\Pb\set{\tX^x(n) > \tX^y(n) \ \text{for at least one value of $n$}} =1.
\end{equation}
\end{lemma}

\begin{proof}
Let $0 < \vae < 1/4$ and $y\le B<\infty$ be such that $H\set{(B,\infty)} = \vae$. Let $w_{\vae}$ be
such that
\[
\abs{\Pb\set{\tX^B(w_{\vae}) \le z} - H(z)} < \vae \quad \text{for all $z$}
\]
and
\[
\abs{\Pb\set{\tX^0(w_{\vae}) \le z} - H(z)} < \vae \quad \text{for all $z$}.
\]
By virtue of Lemma \ref{le3},
\[
\abs{\Pb\set{\tX^x(w_{\vae}) \le z} - H(z)} < \vae \quad \text{for all $z$}.
\]
Write $m$ for the median of the stationary distribution $H$. Obviously,
\begin{align*}
\Pb & \brc{\{B \ge \tX^x(w_\vae) \vee \tX^y(w_\vae)\} \setminus \{B \ge \tX^x(w_\vae) \ge m,
\tX^y(w_\vae) \le m\}}
\\
& \le (1-\vae)^2 - (\tfrac{1}{2}-\vae)^2
\end{align*}
and
\[
(\tfrac{1}{2}-2\vae)^2 < (\tfrac{1}{2}-2\vae)(\tfrac{1}{2}-\vae) \le \Pb \set{B \ge \tX^x(w_\vae)
\ge m, \tX^y(w_\vae) \le m} \le (\tfrac{1}{2}+\vae)^2.
\]

Similarly, for any $j \ge 2$ when $u< v$
\begin{align*}
(\tfrac{1}{2}+\vae)^2 & \ge \Pb \set{\tX^x(j w_\vae) \ge m, \tX^y(j w_\vae) \le m |
\tX^x((j-1)w_\vae)=u, \tX^y((j-1)w_\vae)=v}
\\
& \ge (\tfrac{1}{2}-2\vae)^2
\end{align*}
and
\begin{align*}
\Pb & \brc{\{B \ge \tX^x(jw_\vae) \vee \tX^y(jw_\vae)\} \setminus \{ B \ge \tX^x(jw_\vae) \ge m,
\tX^y(jw_\vae) \le m\}}
\\
& \le (1-\vae)^2 - (\tfrac{1}{2}-\vae)^2 =\frac{3}{4}-\vae.
\end{align*}
Let $T_B = \min\set{j: \tX^x(jw_\vae) \vee \tX^y(jw_\vae) >B}$.

Using previous inequalities, we obtain
\[
\begin{split}
\Pb \set{B \ge \tX^x(jw_\vae) \ge \tX^y(j w_\vae) \; \text{for some $1 \le j < T_B$}} & \ge
\brc{\tfrac{1}{2}-2\vae}^2 \sum_{i=0}^\infty \brc{\tfrac{3}{4}-\vae}^i
\\
& = \frac{\brc{\tfrac{1}{2}-2\vae}^2}{1-\brc{\tfrac{3}{4}-\vae}}
\\
& = \frac{\brc{\tfrac{1}{2}-2\vae}^2}{\tfrac{1}{4}+\vae}.
\end{split}
\]
Letting $\vae\to 0$ completes the proof.
\end{proof}


\begin{lemma} \label{le5} 
Using the same notation as in Lemma \ref{le4},
\[
\Pb \brc{\tX^x(\ell) \ge \tX^y(\ell) \ \text{for some $\ell \le n$}} \xra[n\to\infty]{} 1
\]
uniformly in $0 \le x < y \le B$.
\end{lemma}

\begin{proof}
This follows directly from Lemma \ref{le4} and its proof.
\end{proof}


\begin{lemma} \label{le6} 
Let $\vae > 0$ and let $0< B < \infty$ be such that $H\set{(B,\infty)}<\vae$. Let $B \le A <
\infty$. Then $H_A\set{(B,A)}<\vae$.
\end{lemma}

\begin{proof}
The lemma follows from the fact that $H_A(y) \ge H(y)$ for all $y \ge 0$ (cf.\
Theorem~1 of Pollak and Siegmund, 1986).
\end{proof}


\renewcommand{\proof}{{\sc Proof of Theorem \ref{th1} (i).}}

\begin{proof} 
Let $N_A^{H_A} = \min\set{n: X(n) > A}$ where $X(0) \sim H_A$. By Lemma~\ref{le2},
$N_A^{H_A} \sim \Geom(p_A)$ and
\begin{equation*} 
\lim_{A\to\infty} \Pb\brc{p_A N_A^{H_A} > s} = e^{-s}, \quad s > 0.
\end{equation*}

Let $\vae >0$. Let $0< B < \infty$ be such that $H\set{(B,\infty)}<\vae$. Using the
notation of Lemma~\ref{le4}, let $0< q_B < \infty$ be such that
\begin{equation} \label{bigger1minuseps}
\Pb \brc{\tX^0(n) \ge \tX^B(n) \; \text{for some $n \le q_B$}} > 1- \vae.
\end{equation}

By virtue of Lemma~\ref{le1} and Lemma~\ref{le2}, there exists $A_\vae$ such that
for all $A \ge A_\vae$
\begin{equation} \label{HAH}
\abs{H_A(x)-H(x)} \le \vae \quad \text{for all $0 \le x \le B$}
\end{equation}
and
\begin{equation}  \label{Pexp}
\abs{\Pb\brc{p_A N_A^{H_A} > s} - e^{-s}} \le \vae \quad \text{for all $0 \le s < \infty$}.
\end{equation}

Because the support of $H$ is $[0,\infty)$, it follows from \eqref{HAH} that $p_A q_B
\xra[A\to\infty]{} 0$.

Next, we construct the following sample space. Let $\hX^0(n)$ be a Markov process
(with transition probabilities as $X(n)$) starting at $0$ and let $\hX^B(n)$ be a
Markov process starting at $B$ such that they are independent until the first time
that $\hX^0(n) \ge \hX^B(n)$. Denote this time by $\tau$. After $\tau$, let
$\hX^0$, $\hX^B$ be such that $\hX^0(n) \ge \hX^B(n)$ for all $n \ge \tau$. (This
construction is feasible by virtue of Lemma~\ref{le3} and Lemma~\ref{le4}.)

By virtue of equation \eqref{bigger1minuseps}, $\Pb(\tau \le q_B) > 1-\vae$. Denote
\[
\hN_A^0 = \min\set{n \ge 1: \hX^0(n) > A} \quad \text{and} \quad \hN_A^B =
\min\set{n \ge 1: \hX^B(n) > A}.
\]

Note that $N_A^x$ is stochastically larger than $N_A^y$ if $x < y$.

Now, fix $0 \le s < \infty$ and let $A_B$ be large enough so that $p_A q_B <s$ for
all $A \ge A_B$. Then we have the following chain of equalities and inequalities:
\begin{equation} \label{chain}
\begin{split}
\Pb\brc{p_A N_A^B > s} & = \Pb\brc{p_A \hN_A^B > s}
\\
 & \ge \Pb\brc{p_A \hN_A^B > s, \tau \le q_B}
\\
& = \Pb\brc{p_A \hN_A^B > s, p_A \tau \le p_Aq_B < s}
\\
& = \Pb\brc{p_A \hN_A^B > s, p_A \tau \le p_Aq_B < s, \hN_A^B > \tau}
\\
& \ge \Pb\brc{p_A \hN_A^0 > s, p_A \tau \le p_Aq_B < s, \hN_A^0 > \tau}
\\
& = \Pb\brc{p_A \hN_A^0 > s, p_A \tau \le p_Aq_B < s}
\\
& = \Pb\brc{p_A \hN_A^0 > s, \tau \le q_B}
\\
& \ge \Pb\brc{p_A \hN_A^0 > s} - \Pb\brc{\tau > q_B}
\\
& \ge \Pb\brc{p_A \hN_A^0 > s} - \vae
\\
& = \Pb\brc{p_A N_A^0 > s} - \vae.
\end{split}
\end{equation}
On the other hand,
\begin{equation*}
\begin{split}
\Pb\brc{p_A N_A^B > s} & = \Pb\brc{p_A N_A^{H_A} > s | X(0) =B}
\\
& \le \Pb\brc{p_A N_A^{H_A} > s | X(0)  \le B}
\\
& = \frac{\Pb\brc{p_A N_A^{H_A} > s , X(0)  \le B}}{\Pb\brc{X(0) \le B}}
\\
& = \frac{\Pb\brc{p_A N_A^{H_A} > s , X(0)  \le B}}{H_A([0,B])}
\\
& \le \frac{\Pb\brc{p_A N_A^{H_A} > s}}{H_A([0,B])} .
\end{split}
\end{equation*}
Since by the definition of $B$ and Lemma \ref{le6}, $H_A([0,B]) \ge 1- \vae$, and by equation
\eqref{Pexp}, $P(p_AN_A^{H_A}>s) \le e^{-s} + \vae$, we obtain
\begin{equation} \label{seven}
\Pb\brc{p_A N_A^B > s} \le \frac{e^{-s} + \vae}{1-\vae} .
\end{equation}

Also, since $P(X(0) \ge 0)=H_A([0,A])=1$,
\begin{equation} \label{eight}
\begin{split}
\Pb\brc{p_A N_A^{0} > s} & = \Pb\brc{p_A N_A^{H_A} > s| X(0)=0}
\\
& \ge \Pb\brc{p_A N_A^{H_A} > s| X(0) \ge 0}
\\
& = \frac{\Pb\brc{p_A N_A^{H_A} > s, X(0) \ge 0}}{P(X(0)\ge 0)}
\\
& =\Pb\brc{p_A N_A^{H_A}> s}
\\
& \ge e^{-s} - \vae ,
\end{split}
\end{equation}
where the last inequality follows from equation \eqref{Pexp}.

Putting \eqref{chain} and \eqref{eight} together yields
\begin{equation} \label{nine}
\Pb\brc{p_A N_A^{B} > s} \ge e^{-s} - 2 \vae,
\end{equation}
and putting \eqref{chain} and \eqref{seven} together obtains
\begin{equation} \label{ten}
\Pb \brc{p_A N_A^{0} > s} \le \frac{e^{-s}+\vae}{1-\vae} + \vae .
\end{equation}
Since for all $0 \le x \le B$,
\begin{equation} \label{eleven}
\Pb\brc{p_A N_A^B > s} \le \Pb\brc{p_A N_A^x > s} \le \Pb\brc{p_A N_A^0 > s},
\end{equation}
equations \eqref{nine}--\eqref{eleven} imply that
\[
e^{-s} - 2\vae \le \Pb\brc{p_A N_A^x > s} \le \frac{e^{-s}+\vae}{1-\vae} + \vae \quad \text{for all
$0 \le x \le B$}  .
\]

Finally, fix $x$ and let $\vae\to 0$, so that ultimately $B > x$. This completes
the proof of Theorem~\ref{th1}~(i).
\end{proof}


\

\renewcommand{\proof}{{\sc Proof of Theorem \ref{th1} (ii).}}

\begin{proof} 
Since $N_A^{H_A}$ is distributed $\Geom(p_A)$, $p_A N_A^{H_A}$ has a moment
generating function
\[
M_A^{H_A}(t) = \Eb e^{tp_AN_A^{H_A}}, \quad t < 1,
\]
and it is easy to see that
\begin{equation} \label{Mconv}
M_A^{H_A}(t) \xra[A\to\infty]{}  \frac{1}{1-t} \quad \text{for $t < 1$}.
\end{equation}
Obviously,
\[
M_A^{H_A}(t) = \Eb \Eb \brc{e^{tp_AN_A^{H_A}}| X(0)},
\]
where $X(0)$ has distribution $H_A$. It follows that for every initial state $x \ge 0$ and all $t
<1$ the value of $p_AN_A^x$ has a moment generating function
\[
M_A^x(t) = \Eb e^{t p_A N_A^x}
\]
and
\[
M_A^{H_A}(t) = \Eb M_A^{X(0)}(t) = \int_0^A M_A^x(t) H_A(dx).
\]
For $t \le 0$, by virtue of Theorem \ref{th1}(i)
\[
M_A^x(t) \xra[A\to\infty]{}  \frac{1}{1-t}.
\]

Let $0< \vae < 1$ and $C>0$ be such that $H\{[0,C)\}=\vae$. For fixed $0 < t <1$, let $A(\vae)>C$
be such that
\[
1-\vae < \frac{M_A^{H_A}(t)}{1/(1-t)} < 1+\vae \quad \text{whenever $A \ge A(\vae)$}.
\]
Recall that $X(0)$ has distribution $H_A$, which is a quasi-stationary distribution.

For any $0<\ga <\infty$, Markov's inequality yields
\[
\Pb\brc{M_A^{X(0)}(t) > \ga M_A^{H_A}(t)} \le 1/\ga,
\]
so that for $A \ge A(\vae)$
\begin{equation} \label{twelve}
\Pb\brc{M_A^{X(0)}(t) > \frac{\ga}{1-t}} \le \frac{1+\vae}{\ga} .
\end{equation}
Substituting $\ga=(1+\vae)/\vae$ in \eqref{twelve} yields
\[
\Pb\brc{M_A^{X(0)}(t) > \frac{1+\vae}{\vae}\frac{1}{1-t}} \le \vae .
\]

Since, by Lemma \ref{le6}, $\vae=H\{[0,C)\} \le H_A\{[0,C)\}$, it follows that for $M_A^{X(0)}(t)
\ge \tfrac{1+\vae}{\vae}\tfrac{1}{1-t}$, the value of $X(0)$ cannot exceed $C$. In other words,
\begin{equation} \label{thirteen}
M_A^x(t) \le \frac{1+\vae}{\vae}\frac{1}{1-t} \quad \text{for $x \ge C$ and all $A \ge A(\vae)$} .
\end{equation}

Let $\beta=\min\set{n: X(n) \ge C}$.  Obviously,
\begin{equation} \label{forteen}
M_A^0(t) = \Eb e^{t p_A N_A^0} \le \Eb e^{tp_A\beta} \cdot \Eb e^{tp_AN_A^C} .
\end{equation}

Let $\de_\vae=\Pb\set{X^0(1) \ge C}$. Clearly $\de_\vae \to \Pb\set{X^0(1) > 0} > 0$ as $\vae\to
0$.

Due to the monotonicity of the process $X(n)$, $\beta$ is bounded by a
$\Geom(\de_\vae)$-distributed random variable, so that for $0 < t <1$
\[
1 \le \Eb e^{tp_A\beta} \le \Eb e^{tp_A\Geom(\de_\vae)} = \frac{\de_\vae e^{p_At}}{1- (1-\de_\vae)
e^{p_At}} .
\]
It follows that $\Eb e^{tp_A\beta}$ is bounded as $A\to\infty$ (since $p_A \xra[A\to\infty]{} 0$).
Since $\Eb e^{tp_A N_A^C}=M_A^C(t)$, equations \eqref{thirteen} and \eqref{forteen} imply that
$M_A^0(t)$ is also bounded as $A\to\infty$.

Denote $\varphi(t) = \lim\sup_{A\to\infty} M_A^0(t) < \infty$. Let $\{A_i\}_{i=1}^\infty$ be a
sequence such that $\lim_{i\to\infty}M_{A_i}^0(t)=\varphi(t)$. Construct a set
$\{t_j\}_{j=1}^\infty$ dense in $(0,t)$. Because $M_A^0(u)$ is monotone in $u$, one can obtain a
subsequence $\{A_{ij}\}$ of $\{A_i\}$ such that $M_{A_{ij}}^0(u)$ converges as $j\to\infty$ for all
$0<u<t$. Since the limit is a moment generating function, by Theorem \ref{th1}(i) it must be
$1/(1-t)$. The same argument can be applied to $\lim\inf_{A\to\infty}M_A^0(t)$.

It follows that the limit $\lim_{A\to\infty}M_A^0(t)$ exists and is equal to $1/(1-t)$ for all
$t<1$. Because $M_A^x(t)$ is monotone in $x$ and because of \eqref{Mconv}, $\lim_{A\to\infty}
M_A^x(t)$ necessarily equals $1/(t-1)$ for all $t<1$ and every fixed $x\in [0,\infty)$. This
completes the proof of Theorem \ref{th1}(ii).
\end{proof}

\section{Remarks} \label{s:Remarks}

1. Let $G$ be a distribution with support $[0,A]$ and define the operator $T$ as
\[
T(G) = \text{the distribution of $X(1)$ conditioned on $\set{X(1)\le A, X(0)\sim G}$}.
\]
If $T$ is a continuous operator (in the weak* topology on the distribution functions over $[0,A]$),
then a quasi-stationary distribution exists, i.e., Condition C2 in Theorem \ref{th1} is satisfied
(cf.\ Harris, 1963, Theorem III.10.1).

2. Even if $T$ is not a continuous operator, sometimes Condition C2 can be verified by solving for
$T(G)=G$ and arguing that this is the quasi-stationary distribution. For an example, see Pollak
(1985).

3. The proof can be modified easily to extend Theorem \ref{th1} to the case where the support of
the stationary distribution $H$ is $[c, \infty)$ for some $c>0$ (i.e., the set $[0,c)$ is not in
the state space or is transient).

\section{Examples and Applications} \label{s:Applications}

Theorem \ref{th1} can be applied to a number of popular Harris recurrent Markov
processes. Below we present two examples.  These are of interest when applying
certain change-detection procedures.

\subsection{Example 1: An Additive-Multiplicative Markov Process} \label{ss:S-R}

Let $\La_1, \La_2, \dots$ be non-negative continuous independent and identically distributed
(i.i.d.) random variables with $\beta=\Eb \La_i$ and $\mu = \Eb\log \La_i$. For $x \ge 0$, define
recursively:
\begin{equation} \label{SR}
X(0)=x, \quad X(n) = \brc{1+X(n-1)} \La_n, \quad n=1,2,\dots \; .
\end{equation}
This process is of interest in a number of applications (cf.\ Kesten, 1973; Pollak, 1985, 1987).
For example, in the problem of detecting a change in distribution, the Shiryaev-Roberts statistic
can be written as (cf.\ Pollak, 1985, 1987)
\begin{equation} \label{Rnnew}
R(n)= (1+R(n-1)) \frac{f_{\theta_1}(Y_n)}{f_{\theta_0}(Y_n)}, \quad R(0)=0,
\end{equation}
where $\{Y_n, n \ge 1\}$ are independent, having probability density $f_{\theta_0}$ before a change
and putative density $f_{\theta_1}$ after a change; $\theta_0$ and $\theta_1$ are fixed parameters,
and one stops and declares that the change is in effect at $N_A=\min\{n: R(n) >A\}$.

When $\mu <0$, the process $\{X(n)\}$ is Harris-recurrent and has a stationary distribution (for
any $x\ge 0$). To see this, note that $X(n)$ can be written as
\[
X(n) =\sum_{k=0}^n \prod_{i=k}^n \La_i = \sum_{k=0}^n \exp\set{\sum_{i=k}^n \log \La_i},
\]
where $\La_0=x$. Obviously,
\[
\sum_{k=0}^n \exp\set{\sum_{i=k}^n \log \La_i} \overset{dist}{=} \sum_{k=1}^n \exp\set{\sum_{i=1}^k
\log \La_i} + x \exp\set{\sum_{i=1}^n \log \La_i},
\]
where the right hand-side converges (for every $x \ge 0$ as $n\to\infty$) to the random variable
\[
\sum_{k=1}^\infty \exp\set{\sum_{i=1}^k\log \La_i},
\]
which is a.s.\ finite when $\mu<0$. Since we assumed above that $\La_1$ is continuous, the
quasi-stationary distribution exists (see Remark 1 in Section \ref{s:Remarks}). It follows from
Theorem \ref{th1} that a suitably standardized version of the first exceedance time over $A$ (i.e.,
$p_A N_A^x$) is asymptotically exponentially distributed.

Note that while using the conventional regeneration argument is perhaps possible,
embedding the Markov process \eqref{SR} into ``regenerative cycles" by no means is
either straightforward or obvious, which is especially true when $1\le \beta=\Eb
\La_i < \infty$ and $\mu = \Eb\log \La_i <0$. This case does have meaning for
applications. For example, regard the aforementioned change detection problem. When
there never is a change, the observations $Y_i, i \ge 1$ have density
$f_{\theta_0}$, so that $\beta = \int [f_{\theta_1}(y)/f_{\theta_0}(y)]
f_{\theta_0}(y) dy=1$ while by Jensen's inequality $\mu=\int \log
[f_{\theta_1}(y)/f_{\theta_0}(y)] f_{\theta_0}(y) dy <0$. If there is a change --
for argument's sake let it be in effect from the very beginning -- the observations
$Y_i, i \ge 1$ have density $f_{\theta}$ (not necessarily $f_{\theta_1}$; the
post-change parameter is seldom known in advance, and the putative $\theta_1$ is
merely a representation of a ``meaningful" change). For $\theta$ close to
$\theta_0$, one would obtain $\beta=\int
[f_{\theta_1}(y)/f_{\theta_0}(y)]f_{\theta}(y) dy>1$ and $\mu=\int \log
[f_{\theta_1}(y)/f_{\theta_0}(y)]f_{\theta}(y) dy <0$.

Before going into further details, we discuss an issue related to computing $p_A$, the
standardizing factor. If $p_A$ were amenable to direct calculation, one could use this to
approximate $\Eb N_A^x \approx 1/p_A$. Unfortunately, in most cases direct evaluation of $p_A$ is
not tractable, and evaluation of $\Eb N_A^x$ has to be done by other methods. (But see Pollak,
1985, and Mevorach and Pollak, 1991 for examples that allow some tractability.) \ Nonetheless,
evaluation of $p_A$ is of interest on its own merits (cf.\ Tartakovsky, 2005), as $p_A$ is an
approximation of the probability that there will be a first upcrossing of the threshold $A$ at a
specified time $n$, and $1-(1-p_A)^m$ is an approximation of the probability that there will be a
first upcrossing of $A$ in a given stretch of $m$ observations (i.e., for the ``local false alarm
probability" $\Pb (n \le N_A^x \le n+m-1 | N_A^x \ge n)$). Therefore, if $\Eb N_A^x$ can be
evaluated, $p_A$ can be approximated by $1/\Eb N_A^x$.

Suppose now that $\beta=\Eb \La_i=1$. Let $f_0$ be the density of $\La_i$ and define $f_1(\La)= \La
f_0(\La)$. (Since $\Eb \La=1$, it follows that $f_1$ is a bona fide probability density.) \ Note
that $\La$ is a likelihood ratio, $\La=f_1(\La)/f_0(\La)$. It follows from Pollak (1987) (see also
Tartakovsky and Veeravalli, 2005) that
\begin{equation} \label{ARL-SR}
\Eb_{f_0} N_A^x = \ga^{-1} A (1+o(1)) \quad \text{as $A\to\infty$},
\end{equation}
where $\Eb_{f_0}$ is the expectation with respect to the density $f_0$ and $\ga$ is a constant that
can be calculated by renewal theory (cf.\ Woodroofe, 1982; Siegmund, 1985), so that $p_A\approx
\ga/A$. See Remark in the end of Section \ref{ss:CUSUM} for evaluation of $p_A$ when $\Eb \La_i
\neq 1$.

\subsection{Example 2: A Reflected Random Walk} \label{ss:CUSUM}

Let $\{Z_n\}_{n=1}^{\infty}$ be a sequence of i.i.d.\ continuous random variables with a negative
mean $\mu=\Eb Z_n <0$. For $n \ge 1$, define
\begin{equation} \label{CUSUM}
X(n)= \max\set{0, X(n-1) +Z_n}, \quad X(0)=x \ge 0.
\end{equation}
Since $\mu<0$, the Markov process $\{X(n)\}$ is Harris-recurrent and has a stationary distribution.
To see this, note that
\[
X(n) = \max \set{0,Z_1+\cdots+Z_n +x, Z_2+\cdots+Z_{n-1},\dots,Z_n}.
\]
Write $S_i=\sum_{k=1}^i Z_k$, $S_0=0$. Since the vector $(Z_1,\dots,Z_n)$ has the same distribution
as $(Z_n,\dots,Z_1)$, it follows that
\[
X(n)\overset{dist}{=} \max\set{\max\{0, S_1, S_2,\dots,S_{n-1}\}, x+S_n},
\]
where the right hand-side converges (as $n\to\infty$ for any $x \ge 0$) to the random variable
$\max_{i \ge 0} S_i$, which is a.s.\ finite whenever $\mu=\Eb Z_i <0$.

The process \eqref{CUSUM} describes a broad class of single-channel queuing systems (see, e.g.,
Borovkov, 1976) as well as a popular cumulative sum decision statistic for detecting a change in
distribution (Page, 1954) and has been studied extensively by itself, outside the framework of
general Markov processes. For instance, for $x=0$, the asymptotic exponentiality of the stopping
time
\begin{equation} \label{CUSUMst}
N_a^{x=0} = \min \set{n \ge 1: X(n) >a}, \quad a >0
\end{equation}
(as $a\to\infty$) has been proven by Khan (1995), which can be generalized easily for any $x >0$.
(The process $\{X(n)\}$ obviously is a renewal process, so, although our Theorem~\ref{th1} covers
this example when the conditions C1 and C2 are satisfied, it is not needed to prove asymptotic
exponentiality of $N_A^x$, as it can be derived from general results; cf.\ Asmussen, 2003, Ch.~VI.)

Assume for simplicity that $x=0$. If there exists a positive $\omega$ such that $\Eb e^{\omega
Z_i}=1$, let $f_0(z)$ be the density of $Z_i$ and define $f_1(z)= e^{\omega z} f_0(z)$. Since $\Eb
e^{\omega Z_i}=1$, it follows that $f_1$ is a bona fide probability density, and
$f_1(Z)/f_0(Z)=e^{\omega Z}$ is a likelihood ratio. Hence, assuming that $\mu_1=\int \log
[f_1(z)/f_0(z)] f_1(z) dz < \infty$ and letting
\[
N_a^0 = \min \set{n \ge 1: \max\brc{0, \omega X(n-1) + \omega Z_n} > \omega a},
\]
standard renewal-theoretic methods (cf.\ Woodroofe, 1982; Siegmund, 1985) readily apply to obtain
that
\begin{equation} \label{ARL-CUSUM}
\Eb N_a^0 = \de^{-1} e^{\omega a} (1+o(1)) \quad \text{as $a\to\infty$},
\end{equation}
so that $p_A \approx \de e^{-\omega a}$. Here $0<\de <1$ is a constant that can be computed
explicitly by a renewal-theoretic argument (cf.\ Tartakovsky, 2005).

\begin{remark}
Clearly, $N_a^x$ of Example 2 is larger than $N_A^x$ of Example 1 (with $A=e^a$), so that $\Eb
N_A^x \le \de^{-1} A^\omega (1+o(1))$. Theorem 5 of Kesten (1973) as well as Theorem 4 of Borovkov
and Korshunov (1996) imply that
\[
\Pb (X(\infty) >y) = C/y^{\omega} (1+o(1)) \quad \text{as $y \to\infty$},
\]
where $X(\infty)$ is a random variable that has the stationary distribution of $\{X(n)\}$ and $C$
is a positive finite constant. Note that $X(\infty)$ is stochastically larger than a random
variable that has the quasi-stationary distribution. Therefore, the first upcrossing over $A$ of
the process $X(n)$ starting at a random $X(0)$ distributed like $X(\infty)$ will occur no later
than the first upcrossing over $A$ of the process $X(n)$ starting at a random $X(0)$ that has the
quasi-stationary distribution. The proportion of times that the former exceeds $A$ is
$\Pb(X(\infty) >A)$. It follows that $\Eb N_A^x \ge C^{-1}A^{\omega}(1+o(1))$, so that $p_A$ has an
order of magnitude $1/A^\omega$.
\end{remark}

\subsection{Applications to Sequential Change-Point Detection and a Monte Carlo Experiment} \label{ss:CPD}

The importance of the asymptotic exponentiality of the run length in sequential
change-point detection methods is twofold.  First, it shows that the mean time to
false alarm (the so-called average run length), which is a popular measure of the
false alarm rate, is indeed an exhaustive performance metric. Second, the result
can be used for the evaluation of the local false alarm probabilities of the
corresponding detection schemes (see Example 1 above; see Tartakovsky (2005) for a
more detailed discussion of the importance of local false alarm probabilities in a
variety of applications).

To be more specific, assume that there is a sequence i.i.d.\ variables
(observations) $Y_1, Y_2, \dots$ that follow the density $f_0(y)$ under the
no-change hypothesis (the in-control mode) and the density $f_1(y)$ after the
change occurs (the out-of-control mode). The change occurs at an unknown point in
time $\nu$; $1 \le \nu <\infty$. Therefore, conditioned on $\nu=k$, $Y_n \sim
f_0(y)$ for $n < k$ and $Y_n\sim f_1(y)$ for $n \ge k$. We write $\Pb_\infty$
($\Eb_\infty$) and $\Pb_k$ ($\Eb_k$) respectively for probability measures
(expectations) when there is no change (i.e., $\nu=\infty$) and when the change
occurs at point $k$. Let $Z_n=\log [f_1(Y_n)/f_0(Y_n)]$ be the corresponding
log-likelihood ratio and let $S_n=\sum_{i=1}^n Z_i$. Let $I_1=\Eb_1 Z_1$ and $I_0 =
\Eb_\infty (-Z_1)$ be the Kullback-Leibler information numbers, which are assumed
finite.

We begin with the cumulative sum (CUSUM) test. The CUSUM statistic is given by the recursion
\eqref{CUSUM} and the corresponding stopping time is defined in \eqref{CUSUMst}. The difference
from the previous section is that $Z_n$, $n=1,2,\dots$ are not arbitrary random variables with
negative mean, but rather log-likelihood ratios with mean $\mu=-I_0$. This simplifies most of the
calculations, since $\Eb e^{Z_n}=1$. Recall that in this section we denote this expectation by
$\Eb_\infty$.

Rewrite the corresponding stopping time in the following form
\begin{equation} \label{Wst}
N_A = \min \set{n \ge 1: \max\set{1, W(n-1)+ e^{Z_n}} >  A} ,
\end{equation}
where $W(0)=1$ and $A=e^{a}$. The asymptotic approximation for the average run length to false
alarm \eqref{ARL-CUSUM} holds with $\omega=1$, $e^a=A$, and $\de=I_1 \ga^2$ (cf.\ Tartakovsky,
2005), which implies that $p_A \sim I_1 \ga^2/A$. Here $\ga = \lim_{y\to\infty} \Eb_1
\exp\{-(S_{\tau_y}-y)\}$, where $\tau_y=\min\{n: S_n > y\}$ is the first time when the random walk
$S_n=\sum_{i=1}^n Z_i$ crosses the level $y$. The constant $\ga$ is the subject of renewal theory
(cf.\ Woodroofe, 1982 or Siegmund, 1985) and can be computed explicitly.

We now proceed with the Shiryaev-Roberts detection test. The Shiryaev-Roberts statistic is defined
by \eqref{Rnnew}, where $\tfrac{f_{\theta_1}(Y_n)}{f_{\theta_0}(Y_n)}=e^{Z_n}$ and $R(0)=0$. The
corresponding stopping time is
\[
\hN_A = \min \set{n \ge 1: R(n) > A}.
\]
We now denote it by $\hN_A$ to distinguish from the CUSUM stopping time in the following
calculations and comparison.

Since $\Eb_\infty e^{Z_n}=1$, the process $R(n)-n$ is a zero-mean martingale, which allows us to
approximate the average run length to false alarm:
\[
\Eb_\infty \hN_A \sim \ga^{-1} A \quad \text{as $A \to \infty$}.
\]
This approximation follows from \eqref{ARL-SR} above. The distribution of the
Shiryaev-Roberts stopping time is approximately $\Exp(p_A)$ with $p_A \sim \ga/A$.
(The asymptotic exponentiality of the suitably standardized run length to false
alarm has been shown by Yakir, 1995.)

In order to verify the accuracy of asymptotic approximations for reasonable values of the threshold
$A$, we performed Monte Carlo (MC) simulations for the following example. Consider the case where
observations are independent, originally having an $\Exp(1)$ distribution, changing at an unknown
time to $\Exp(1/(1+q))$, i.e.,
\begin{equation} \label{Exppmodel}
f_0(y)= e^{-y}\Ind{y\ge 0}, \quad f_1(y) = \frac{1}{1+q}e^{-y/(1+q)}\Ind{y\ge 0}, \quad q >0.
\end{equation}
In this case
\[
I_1=q-\log (1+q) \quad \text{and} \quad \ga =1/(1+q).
\]
Applying Example 1, the likelihood ratio is $\La_n= e^{Z_n} =(1+q)^{-1}e^{q Y_n/(1+q)}$ and the
average run length (ARL) to false alarm of the Shiryaev-Roberts procedure is
\begin{equation} \label{ARL-SRexp}
\ARL_{{\mrm SR}}(A)=\Eb_\infty \hN_A \approx (1+q) A.
\end{equation}

Applying Example 2, an approximation of the ARL to false alarm of the CUSUM test is
\begin{equation} \label{ARL-CUSUMexp}
\ARL_{{\mrm CU}}(A)=\Eb_\infty N_A \approx \frac{(1+q)^2}{q-\log(1+q)} A .
\end{equation}

\begin{table}[ht!]
\begin{center}
\caption{The ARL versus threshold for the CUSUM test for $q=3$}
\label{t:ARLforCUSUM} {\small
\begin{tabular} {|c||*{9}{c|}} 
\hline\hline $A$& 1.2&1.7&2.5&4.6&9.2&13.0&17.1&21&41
\\ \hline
FO $\ARL_{{\mrm CU}}$& 11.90&16.86&24.79&45.61&91.22&128.90&169.55&208.22&406.52
\\ \hline
HO $\ARL_{{\mrm CU}}$& 7.96&12.36&19.69&39.56&84.07&121.21&161.43&199.77&397.02
\\ \hline
MC $\widehat{\ARL}_{{\mrm CU}}$&
8.04&12.45&19.79&39.57&84.33&121.23&161.88&200.44&397.16
\\ \hline
MC $\SD(N_A)$& 7.49&11.88&19.18&38.61&83.21&119.73&159.91&198.97&396.84
\\ \hline\hline
\end{tabular}}
\end{center}
\end{table}

\begin{table}[ht!]
\begin{center}
\caption{The ARL versus threshold for the Shiryaev-Roberts test for $q=3$}
\label{t:ARLforSR} {\small
\begin{tabular} {|c||*{9}{c|}} 
\hline\hline $A$& 1&2&5&10&20&30&40&50&100
\\ \hline
$\ARL_{{\mrm SR}}$& 4&8&20&40&80&120&160&200&400
\\ \hline
MC $\widehat{\ARL}_{{\mrm SR}}$&
4.01&8.03&20.00&39.94&79.99&119.82&159.17&200.42&399.46
\\ \hline
MC $\SD(\hN_A)$& 3.00&6.78&18.34&37.92&77.33&117.39&157.19&197.90&396.94
\\ \hline\hline
\end{tabular} }
\end{center}
\end{table}

We simulated the CUSUM and Shiryaev-Roberts procedures under the assumption of no change (i.e., all
simulated observations are $\Exp(1)$). Each combination of (test,threshold) was simulated 100,000
times. The results are reported in Tables \ref{t:ARLforCUSUM} and \ref{t:ARLforSR}. We present the
results of simulations when the parameter $q=3$, which is a reasonable value in certain
applications such as detection of a randomly appearing target in noisy measurements, in which case
$q$ is the signal-to-noise ratio (see, e.g., Tartakovsky, 1991 and Tartakovsky and Ivanova, 1992).
It is seen that the approximation \eqref{ARL-SRexp} for the Shiryaev-Roberts test is very accurate
for all threshold values, even when the ARL is small. On the other hand, the approximation
\eqref{ARL-CUSUMexp} for the CUSUM test (given in the row ``FO $\ARL_{{\mrm CU}}"$ in Table
\ref{t:ARLforCUSUM}, where FO stands for ``first order") is not especially accurate. This happens
primarily because the first order approximation takes into account only the first term of expansion
and ignores the second term $O(\log A)$ as well as constants. An accurate, higher order (HO)
approximation can be obtained using the results of Tartakovsky and Ivanova (1992) which give:
\begin{align*}
\ARL_{{\mrm CU}}(A) & \approx \frac{(1+q)^2}{q- \log (1+q)} A - \frac{1}{\log (1+q)-q/(1+q)} \log A
\\
& \quad - \frac{1+q}{q- \log (1+q)} - \frac{q}{(1+q)\log(1+q)-q}.
\end{align*}

In Table \ref{t:ARLforCUSUM}, the row ``HO $\ARL_{{\mrm CU}}$" corresponds to this
latter approximation, which perfectly fits the MC estimates (denoted by ``MC
$\widehat{\ARL}_{{\mrm CU}}$") for all tested threshold values $A\ge 1.2$.

In these tables we also present the MC estimates of standard deviations $\SD(N_A)$
and $\SD(\hN_A)$ of the stopping times. As one would expect, the standard
deviations are the same (approximately) as the means, and the similarity grows as
$A$ increases. The fit is slightly better for the CUSUM test.

\begin{figure}[!ht]
 \centering
   \subfigure[CUSUM test: $q=3$, $A=13$]{
  \includegraphics[width=0.46\textwidth]{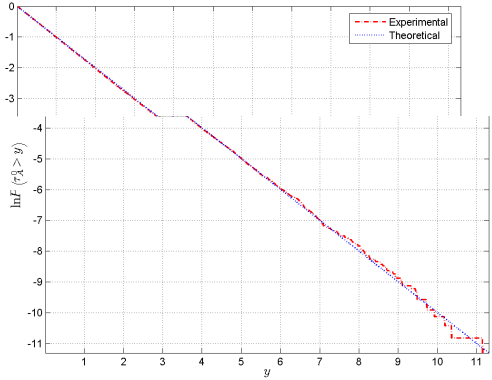}
  \label{f:logCDFforCUSUM}
  }
  \subfigure[Shiryaev-Roberts test: $q=3$, $A=40$]{
  \includegraphics[width=0.46\textwidth]{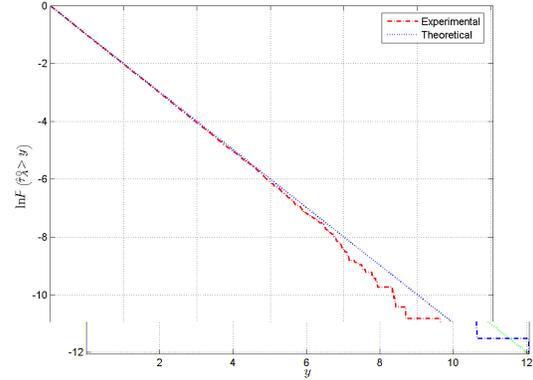}
\label{f:logCDFforSR}
  }
 \caption{Empirical estimates of $\log[\Pb_\infty(\tau_A > y)]$ and $\log[\Pb_\infty(\hat{\tau}_A > y)]$ for
 the CUSUM and Shiryaev-Roberts procedures}
  \label{f:logCDF}
\end{figure}

\begin{figure}[!ht]
\centering
  \subfigure[CUSUM test: $q=3$, $A=13$]{
  \includegraphics[width=0.46\textwidth]{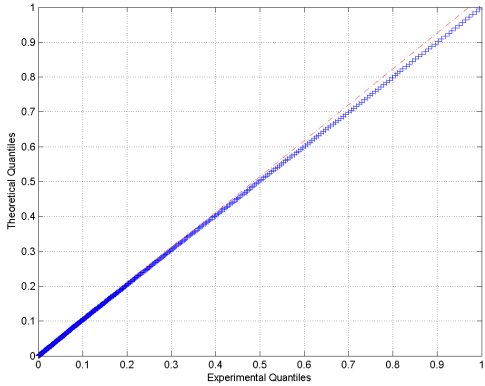}
\label{f:QQplotforCUSUM}
 }
 \subfigure[Shiryaev-Roberts test: $q=3$, $A=40$]{
  \includegraphics[width=0.46\textwidth]{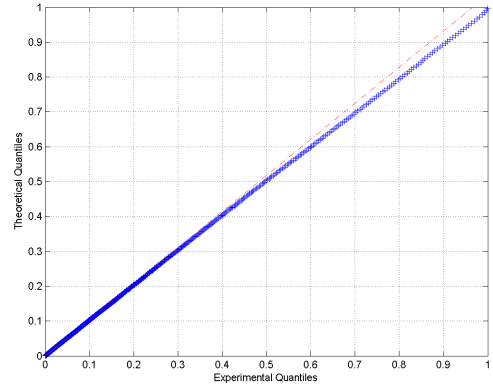}
  \label{f:QQplotforSR}
 }
 \caption{QQ-plots for the stopping times of the CUSUM and Shiryaev-Roberts procedures}
 \label{f:QQplots}
\end{figure}

Figures \ref{f:logCDFforCUSUM} and \ref{f:logCDFforSR} show the logarithm of the
empirical (MC estimates) survival functions $\log \Pb_\infty(\tau_A>y)$ and $\log
\Pb_\infty(\hat{\tau}_A>y)$ for the CUSUM and Shiryaev-Roberts procedures, where
$\tau_A=N_A/\widehat{\ARL}_{{\mrm{CU}}}$ and
$\hat{\tau}_A=\hN_A/\widehat{\ARL}_{{\mrm{SR}}}$ are the corresponding standardized
stopping times, along with the logarithm of the exponential probability plot $\log
e^{-y}=-y$. The quantile-quantile plots (QQ-plots) for the stopping times are shown
in Figures \ref{f:QQplotforCUSUM} and \ref{f:QQplotforSR}. The QQ-plots display
sample quantiles of $N_A$ and $\hN_A$ versus theoretical quantiles from the
exponential distribution. If the distributions of the stopping times are
exponential, the plots will be close to linear. These figures show that, for the
chosen putative value of the post-change parameter ($q=3$), the exponential
distribution approximates the distributions of the stopping times very well. It is
seen that the exponential approximation works very well already for $A=13$
($\ARL_{{\mrm CU}} \approx 120$) for the CUSUM test and for $A=40$ ($\ARL_{{\mrm
SR}} \approx 160$) for the Shiryaev-Roberts test. When considering that in
practical applications the values of the ARL to false alarm usually range from 300
and upwards, the exponential distribution seems to be a perfect fit.

\section*{Acknowledgements}

We are grateful to Alexei Polunchenko for the help with Monte Carlo simulations.

Moshe Pollak is Marcy Bogen Professor of Statistics at the Hebrew University of
Jerusalem. His work was supported in part by a grant from the Israel Science
Foundation, by the Marcy Bogen Chair of Statistics at the Hebrew University of
Jerusalem, and by the U.S.\ Army Research Office MURI grant W911NF-06-1-0094 at the
University of Southern California. The work of Alexander Tartakovsky was supported
in part by the Marcy Bogen Chair of Statistics at the Hebrew University of
Jerusalem and by the U.S.\ Office of Naval Research grant N00014-06-1-0110 and the
U.S.\ Army Research Office MURI grant W911NF-06-1-0094 at the University of
Southern California.

\end{document}